\documentclass[a4paper,10pt]{article}

\usepackage[utf8]{inputenc}
\usepackage{amsmath, amsfonts, amssymb, amsthm} \parskip1pt

\usepackage{mathrsfs}
\usepackage{amsfonts,graphics,epsfig}
\usepackage{amsmath, amsfonts, amssymb, amsthm, amscd}

\newcommand{\SL}{\operatorname{SL}}
\newcommand{\GL}{\operatorname{GL}}
\newcommand{\so}{\operatorname{so}}
\newtheorem{theorem}{Theorem}[section]

\newtheorem{remark}[theorem]{Remark}

\newtheorem{corollary}[theorem]{Corollary}

\newtheorem{acknowledgement}{Acknowledgement}

\def\R{\mathbb R} \def\Z{\mathbb Z}

\def\Q{\mathbb Q}

\def\Aut{\rm Aut}

\date{}

\def\build#1_#2^#3{\mathrel{\mathop{\kern 0pt#1}\limits_{#2}^{#3}}}

\let\ap=\alpha

\def\smallsquare{\vbox{\hrule\hbox{\vrule height 1 ex\kern 1
ex\vrule}\hrule}}

\title{Nilpotent Lie groups and hyperbolic automorphisms}
\author{Manoj Choudhuri and C. R. E. Raja}

\date{}

\begin{document}
\maketitle
\begin{abstract}
 A connected Lie group admitting an expansive automorphism is known to be nilotent, but all nilpotent Lie groups do not admit expansive automorphism. In this article, we find sufficient conditions for a class of nilpotent Lie groups to admit expansive automorphism.
\end{abstract}

Let $G$ be a locally compact (second countable) group and $\Aut (G)$
be the group of (continuous) automorphisms of $G$.  We consider
automorphisms $\ap$ of $G$ that have a specific property that there
is a neighborhood $U$ of identity such that $\cap _{n \in \Z} \ap ^n
(U) = \{ e \}$ - such automorphisms are known as expansive.
The structure of compact connected groups admitting expansive automorphism has been studied in \cite{Ei1}, \cite{Ei2}, \cite{Wu}, \cite{Lam} and it was shown in \cite{Lam} that a compact connected group admitting an expansive automorphism is abelian. The same problen has been studied for totally disconnected groups as well, see \cite{Ki} and \cite{Wi} for compact totally disconnected groups and a recent paper of Helge Glockner and the second named author (\cite{GlR}) for locally compact totally disconnected groups.
In \cite{Ri}, Riddhi Shah has shown that a connected locally compact group admitting expansive automorphism is nilpotent. The result was known previously for Lie groups (see Proposition $7.1$ of \cite{GlR} or Proposition $2.2$ of \cite{Ri} for a proof of this fact and further references). But all nilpotent Lie groups do not admit expansive automorphism which can be seen easily by considering the circle group. In this article, we try to classify a class of nilpotent Lie groups which admit expansive automorphism. Since in the Lie theory context, expansive automorphisms are better known as hyperbolic automorphisms, from now on we also use the term hyperbolic automorphisms instead of expansive automorphisms.

Let $N$ be $k$-step nilpotent Lie algebra generated by $q$ linearly independent elements $e_1,...,e_q$ (see \cite{Sato} for details about structure of nilpotent Lie algebras). If $U$ is the vector space generated by $e_1,...,e_q$, then $N=U\oplus [N,N]$. We denote by $N^i$ the $i$th term in the central series of $N$, that is
$$N^i=[N,N^{i-1}].$$ Let $\widetilde{N}$ denote the free $k$-step nilpotent Lie algebra on $q$ generators. Then there exists a
Lie algebra homomorphism $T$ (say) from $\widetilde{N}$ to $N$ such that $T(e_i)=e_i$. Let $\text{Ker}(T)=\widetilde{A}$. We say
an element $[e_{i_1},[e_{i_2},[...[e_{i_{m-1}},e_{i_m}]...]$ a monomial of lenght $m$. $e_i$ is a monomial of lenght $1$.
$\widetilde{N}$ is generated by the monomials of lenght less or equal to $k$. We call these monomials the standard basis elements
of $\widetilde{N}$.

Now let $\alpha\in \SL(q,\R)$. Then by defining
$$\alpha([e_{i_1},[e_{i_2},[...[e_{i_{m-1}},e_{i_m}]...])=[\alpha(e_{i_1}),[\alpha(e_{i_2}),[...[\alpha(e_{i_{m-1}}),
\alpha(e_{i_m})]...],$$
$\alpha$ can be extended to an automorphism of $\widetilde{N}$, and we denote this automorphism of $\widetilde{N}$ by the same $\alpha$. Now $\alpha\in\SL(q,\R)$ is an automorphism of $N$ if $\alpha(\widetilde{A})\subset \widetilde{A}$ and $$\overline{G}=\{\alpha\in\SL(q,\R)|\alpha(\widetilde{A})\subset \widetilde{A}\}$$ is an algebraic subgroup of $\SL(q,\R)$. Moreover, if the basis elements of $\widetilde{A}$ are rational linear combinations of the monomials (the standard basis elements of $\widetilde{N}$), then $\overline{G}$ is defined over $\Q$. This can be seen by completing a basis of $\widetilde{A}$ to a basis of $\widetilde{N}$ and then looking at the representation of $G$ in $\GL(\widetilde{N})$ with respect to this basis.

Now if $\alpha\in \overline{G}$, then $\alpha$ is an automorphism of $N$. It is easy to see that $\alpha$ is an automorphism of $N^k$ as well. This gives rise to a representation of $\overline{G}$ in $\GL(N^k)$. This representation is defined over $\Q$, which can be seen by first looking at the representation of $\overline{G}$ in $\GL(N)$ and then restricting it to $N^k$.


\begin{theorem}\label{main}
Let $N$, $\widetilde{N}$, $\widetilde{A}$, $T$ be as above such that $k<q$. Let $\Gamma=\Z w_1\oplus ... \oplus \Z w_d $ be a vector space lattice in $N^k$ such that $w_1,...,w_d$ are rational linear combinations of a $\Q$-basis of $N^k$. If $G$ is an algebraic subgroup of $\overline{G}$ defined over rationals such that $G$ is semi-simple and all the weights of $G$ for its action on $N$ are non-trivial, then there is a hyperbolic automorphism $\alpha$ of $N$ such that $\alpha$ preserves the lattice $\Gamma$.
\end{theorem}
\begin{proof}

Sice $G$ is semisimple and defined over $\Q$, it has no non-trivial character defined over $\Q$. Then it follows from Theorem $9.4$ of \cite{BH} that the set of $\Z$-points of $G$, denoted by $G(\Z)$ is a lattice (i.e., a discrete subgroup of $G$ such that $G/G(\Z)$ carries a $G$-invarient finite measure) in $G$. Note that if we rewrite the representation $\rho$ in terms of the generating set of $\Gamma$, it remains a polynomial representation defined over $\Q$ because of the rationality assumption on $\Gamma$. From now on we are going to consider $\rho$ as this representation only. Now it follows from Proposition $10.13$ of \cite{Rg} that there is a finite index subgroup $G_m(\Z)$ of $G(\Z)$ such that $\rho(G_m(\Z))$ preserves the lattice $\Gamma$.
Our theorem will be proved if we can show that $G_m(\Z)$ contains some hyperbolic automorphism of $N$.

Let $D$ be a maximal $\R$-split torus in $G$ and let $C$ be the centralizer of $D$ in $G$. Then $C$ is a Cartan subgroup of $G$ and it follows from Theorem $2.8$ of \cite{PR1} that up to conjugation by an element of $G$,  $L=C\cap G_m(\Z)$ is a uniform lattice in $C$. Since $G$ is semisimple, without loss of generality, we may assume that $C=KD$, where $K$ is compact, $K$ and $D$ commute. Then if $\pi$ denotes the projection map from $C$ to $D$, $\pi(L)$ is a uniform lattice in $D$. By Borel density theorem (see \cite{BO2} for details), $L$ being a lattice, is Zariski dense in $D$. Let $D^h$ be the set of elements of $D$, which acts as hyperbolic automorphisms on $N$.

Now suppose $\alpha\in D$ and as an element of $\SL(q,\R)$, let $\alpha_1,\ldots,\alpha_q$ be the eigenvalues of $\alpha$. If we consider the action of $\alpha$ as automorphism of $\widetilde{N}$, then it is easy to see that the eigenvalues of $\alpha$ are given by the set $$S=\{\alpha_{i_1}\alpha_{i_2}\ldots \alpha_{i_n}\}_{1\leq n\leq k }.$$ Since $\alpha$ preserves the subspace $\widetilde{A}$ and $N\cong \widetilde{N}/\widetilde{A}$, it follows that there is a subset $S_1$ of $S$ such that  for the action of $\alpha$ on $N$, the eigenvalues are given by the elements of the set $S_1$. Let us denote the elements of $S_1$ by $s_1,s_2,\ldots s_N$. Then $$D^h=\{\alpha\in D|s_i\neq 1,1\leq i\leq N\}.$$ Since we have assumed that all the weights of $G$ (which means all the weights of a maximal $\R$-split torus of $G$) are non-trivial, it follows that $D^h$ is non-empty.
Clearly, $D^h$ is a Zariski open subset of $D$. Then there is some $\alpha\in L\cap D^h$ and the theorem is proved.
\end{proof}

\begin{corollary}\label{cor1}
 Along with all the assumptions of Theorem \ref{main}, let $\Gamma'$ be another lattice in $N^k$ given by $$\Gamma'=\rho(g)\Gamma$$ for some $g\in\overline{G}$. Then there is a hyperbolic automorphism $\alpha'$ of $N$ such that $\alpha'$ preserves $\Gamma'$.
\end{corollary}
\begin{proof}
 By the above theorem, let $\alpha$ be the automorphism of $N$ which preserves $\Gamma$. Let $\alpha'=g\alpha g^{-1}$. Then it is clear that $\alpha'$ is a hyperbolic automorphism of $N$ which preserves $\Gamma'$.
\end{proof}
\begin{remark}
 If $N$ is a free $k$-step nilpotent Lie algebra on $q$ generators, then $\SL(q,\R)\subset \text{Aut}(N)$ and, therefore, for any rationally defined full latice $\Gamma$ in $N^k$ there is a hyperbolic automorphism of $N$ preserving $\Gamma$. Therefore, the corresponding Lie group admits expansive automorphisms. If $N$ is a $2$-step nilpotent Lie algebra then it admits a positive graduation and then it is easy to produce hyperbolic automorphisms (see \cite{S} for details) of $N$, and therefore the simply connected Lie group whose Lie algebra is $N$ admits hyperbolic automorphisms. But all simply connected nilpotent Lie groups do not admit hyperbolic automorphism. In \cite{DL}, an example of a $3$-step nilpotent Lie algebra $\mathfrak{g}$ is discussed, all of whose derivations are nilpotent.  Thus, if $G$ is the simply
 connected nilpotent Lie group with Lie algebra $\mathfrak{g}$,
then the connected component of $\Aut (G)$ is an unipotent group.
Since $\Aut (G)$ is an algebraic group, it has only finitely many
connected components and hence no automorphism of $G$ is hyperbolic.
\end{remark}
Now we discuss an example where the hypothesis of Theorem \ref{main} is verified and which is different from the trivial examples discussed in the above remark. Let $N$ be a $2$-step nilpotent Lie algebra such that $[N,N]$ has dimension $p$ and $U$ has dimension $q$. Following \cite{E1} and \cite{E2}, we call $N$ a $2$-step nilpotent Lie algebra of type $(p,q)$. It was shown in \cite{E1} that any $2$-step nilpotent Lie algebra of type $(p,q)$ is isomorphic to a Lie algebra $\mathfrak{N}=\R^q\oplus W$, which is called a standard metric nilpotent Lie algebra of type $(p,q)$ (see \cite{E1} for details) where $W$ is a $p$-dimensional subspace of $\so(q,\R)$. If $W=\so(q,\R)$, then $\mathfrak{N}$ is the free $2$-step nilpotent Lie algebra on $q$ generators. $\SL(q,\R)$ acts on $\R^q\oplus \so(q,\R)$ by automorphism where the action is given as follows: for $g\in\SL(q,\R)$ and $(v,w)\in \R^q\oplus \so(q,\R)$, $$g(v,w)=(gv,gwg^t).$$ Also $$\SL(q,\R)_W=\{g\in\SL(q,\R):gwg^t\in W, \hspace*{0.1cm} \text{whenever}
\hspace*{0.1cm} w\in W\}$$ acts on $\R^q\oplus W$ by automorphism. We denote $\SL(q,\R)_W$ by $G_W$.

The action of $\SL(q,\R)$ on $\so(q,\R)$ given by
$$ g(w)=gwg^t \eqno (1)$$ for $g\in \SL(q,\R)$ and $w\in \so(q,\R)$, gives rise to the representation $\rho$ of $\SL(q,\R)$ in $\SL(\so(q,\R))$, which is an algebraic representation defined over $\Q$. Let $g\in\SL(q,\R)$ be given by the matrix $(g_{ij})$. Let us try to write the matrix of $\rho(g)$ in the standard orthonormal basis of $\so(q,\R)$. $\rho(g)$ is a $d\times d$ matrix where $d=\frac{1}{2}q(q-1)$. Each row (column) of $\rho(g)$ corresponds to a member of the double indexed set $I=\{ij:1\leq i<j\leq q\}$ which has $d$ number of elements because of the following: if $\{e_1,...,e_q\}$ denotes the standard basis of $\R^q$, then
\begin{align}\label{equality1}
e_{ij}:= [e_i,e_j]
\end{align}
is one of the elements of the standard orthonormal basis of
$\so(q,\R)$ with respect to the inner product described in \cite{E1}, and
each element of the standard orthonormal basis of $\so(q,\R)$ is
determined by \ref{equality1}.  Now let $\rho(g)_{ij}$ be an entry
of the matrix $\rho(g)$. Then $i$ corresponds to two rows of the
matrix $g$ and we denote them by $i_1$ and $i_2$. Also, $j$
corresponds to two columns of the matrix $g$ which we denote by
$j_1$ and $j_2$. Then a straightforward calculation shows that
$\rho(g)_{ij}=g_{i_1j_1}g_{i_2j_2}-g_{i_2j_1}g_{i_1j_2}$. That is
$\rho(g)_{ij}$ is the determinant of the $2\times 2$ minor of $g$ formed by the $i_i$th, $i_2$th rows and $j_1$th, $j_2$th columns of
$g$. We call a $p$-dimensional subspace $W$ to be standard subspace if $W$ is the span of $p$ elements from the standard orthonormal basis of $\so(q,\R)$.

Now let $p=3$ and $q$ be such that $q-3\geq 3$. Let $W$ be the standard $3$-dimensional subspace of $\so(q,\R)$ spanned by the first $3$ elements of the standard orthonormal basis of $\so(q,\R)$ i. e., $W$ is spanned by $e_{12}$, $e_{13}$ and $e_{23}$. Also let $\Gamma=\Z e_{12}\oplus \Z e_{13}\oplus \Z e_{23}$. Then it is easy to see that $G=\SL(3,\R)\times \SL(q-3,\R)$ is a subgroup of $G_W$ satisfying the hypothesis of Theorem \ref{main}. Therefore, there is a hyperbolic automorphism $\alpha$ of $\mathfrak{N}=\R^q\oplus W$ such that $\alpha$ preserves the lattice $\Gamma$. Hence the nilpotent Lie group $\mathfrak{N}/\Gamma=\R^q\oplus W/\Gamma$ admits a hyperbolic automorphism.

\begin{acknowledgement}
The authors would like to thank Prof. S. G. Dani for suggesting
\cite{PR1} and Prof. E. Vinberg for pointing out the mistakes in the
previous version. The first named author acknowledges the support of
National Board for Higher Mathematics, India, through NBHM
post-doctoral fellowship during which the major portion of the work
is done.
\end{acknowledgement}
\begin{center}
{\bf Summary of revisions over the previous version}
\end{center}
Note that some mistakes were found in the previous version of this paper titled\\ `` 2-step nilpotent Lie groups and hyperbolic automorphisms ". The mistakes were found in Lemma $4.2$. It was also found out that we need some extra assumption in Theorem $1.1$. In the present version, we have put this extra assumption in Theorem $0.1$, but we have not restricted ourselves to $2$-step nilpotent Lie groups, rather we have considered $k$-step ($k\geq2$) nilpotent Lie groups. 
 
As Lemma $4.2$ of the previous version fails to hold, in the present version, we do not use the density arguments as we did in the previous version, to achieve same kind of results for a more general class of nilpotent Lie groups apart from those mentioned in Theorem $0.1$ and Corollary $0.2$. Therefore, we think the results in Section $3$ of the previous version regarding convergence in the space of closed subgroups and the use of them in Section $4$ become irrelevant for the present version. So, we have not included those results in the present version, we plan to put those results in a suitable place of a seperate article. 
 
\smallskip
\bibliography{refnil}
\bibliographystyle{plain}

\begin{itemize}
 \item Manoj Choudhuri,\\
 Institute of Infrastructure Technology Research and Management,\\
 Near Khokhara Circle, Maninagar (East), Ahmedabad 380026,\\
 Gujarat, India.\\
 email: manojchoudhuri@iitram.ac.in
 \item C. R. E. Raja,
 Statistics and Mathematics Unit,\\
 Indian Statistical Institute Bangalore Center,\\
 8th Mile, Mysore Road, R. V. C. E. Post,\\
 Bangalore 560059, Karnataka, India.\\
 email: creraja@isibang.ac.in

\end{itemize}

\end{document}